\documentclass[12pt,reqno]{amsart}
\usepackage[T1]{fontenc}
\usepackage{lmodern,amsmath,amssymb,amsthm,mathtools,mathrsfs}
\mathtoolsset{showonlyrefs=true}
\usepackage[letterpaper,margin=1in]{geometry}
\usepackage{microtype,array,enumitem,xcolor}
\usepackage[pdftitle={Weighted composition operators between Qp spaces: boundedness and compactness},pdfsubject={Boundedness and compactness: a research manuscript},colorlinks=true,linkcolor=blue!45!black,citecolor=blue!45!black,urlcolor=blue!45!black]{hyperref}

\allowdisplaybreaks[2]
\makeatletter
\renewcommand\subsection{\@startsection{subsection}{2}%
  \z@{.7\linespacing\@plus.5\linespacing}{.35\linespacing}%
  {\normalfont}}
\makeatother
\newtheorem{theorem}{Theorem}[section]
\newtheorem{proposition}[theorem]{Proposition}
\newtheorem{lemma}[theorem]{Lemma}
\newtheorem{corollary}[theorem]{Corollary}
\theoremstyle{definition}

\theoremstyle{remark}
\newcommand{\D}{\mathbb D}\newcommand{\T}{\mathbb T}\newcommand{\C}{\mathbb C}

\newcommand{\ind}{\mathbf 1}

\newcommand{\Hol}{\operatorname{Hol}}
\newcommand{\supp}{\operatorname{supp}}\newcommand{\id}{\operatorname{id}}

\begin{document}

\title[Weighted composition operators between $Q_p$ spaces]{Weighted composition operators between $Q_p$ spaces: boundedness and compactness}

\author{Wujun Cao, Zhouyuan Jiang and Songxiao Li}
\address{Wujun Cao: Department of Mathematics, Shantou University, 515063, Shantou, Guangdong, P.R. China.}
\email{24wjcao@stu.edu.cn}
\address{Zhouyuan Jiang: Department of Mathematics, Shantou University, 515063, Shantou, Guangdong, P.R. China.}
\email{25zyjiang@stu.edu.cn}
\address{Songxiao Li: Department of Mathematics, Shantou University, 515063, Shantou, Guangdong, P.R. China.}
\email{jyulsx@163.com}

\date{}

\begin{abstract}
We characterize boundedness and compactness of weighted composition
operators $W_{u,\varphi}:Q_{p_1}\to Q_{p_2}$, where
$W_{u,\varphi}f=u(f\circ\varphi)$, $p_1>0$, and $0<p_2<1$.
A separation theorem for the two terms in $(W_{u,\varphi}f)'$,
together with value and derivative embedding theorems, gives two
capacity conditions on $u$ and $\varphi$. For $p_1=1$ and $p_1>1$,
the derivative condition has a maximal-integral and an area-integral
form  respectively.
\end{abstract}
\maketitle

\section{Introduction}\label{sec:intro}

Let $\D=\{z\in\C:|z|<1\}$, let $\T=\partial\D$, and let
$dA(z)=dx\,dy/\pi$ be normalized area measure. Write $\Hol(\D)$
for the analytic functions on $\D$.
For $a\in\D$, set
\begin{equation}\label{eq:auto}
 \varphi_a(z)=\frac{a-z}{1-\overline a z},\qquad
 1-|\varphi_a(z)|^2
 =\frac{(1-|a|^2)(1-|z|^2)}{|1-\overline a z|^2}.
\end{equation}
For $p>0$, the space $Q_p$ consists of all $f\in\Hol(\D)$ such that
\begin{align}
 \|f\|_{Q_p}&=|f(0)|+\|f\|_{Q_p,*}<\infty,\label{eq:Qnorm}\\
 \|f\|_{Q_p,*}^2
 &=\sup_{a\in\D}\int_\D|f'(z)|^2
                 (1-|\varphi_a(z)|^2)^p\,dA(z).
 \label{eq:Qseminorm}
\end{align}
We use this norm throughout. It is equivalent to the usual norm
defined by the Green function $\log(1/|\varphi_a(z)|)$; see
\cite[Section 2.2]{WZ}. The space $Q_p$ is a Banach space and its
seminorm is M\"obius invariant: for every automorphism $\psi$ of $\D$,
\[
 \|f\circ\psi\|_{Q_p,*}=\|f\|_{Q_p,*}.
\]
The scale contains two classical spaces:
\begin{equation}\label{eq:classicalspaces}
 Q_1=\mathrm{BMOA},\qquad Q_p=\mathcal B\quad(p>1),
\end{equation}
with equivalent norms, where BMOA is the analytic subspace of
boundary BMO and $\mathcal B$ is the Bloch space with norm
\[
 \|f\|_{\mathcal B}=|f(0)|+
                \sup_{z\in\D}(1-|z|^2)|f'(z)|.
\]
Moreover $Q_r\subseteq Q_p$ when $0<r\leq p$.
We write $\mathring Q_p=\{f\in Q_p:f(0)=0\}$.

Let $u\in\Hol(\D)$ and let $\varphi:\D\to\D$ be holomorphic. The
weighted composition operator induced by $u$ and $\varphi$ is
\begin{equation}\label{eq:Wdef}
 W_{u,\varphi}f(z)=u(z)f(\varphi(z)).
\end{equation}
For $u=1$ it is the composition operator $C_\varphi$, and for
$\varphi=\id$ the multiplication operator $M_u$. We seek necessary
and sufficient conditions on $u$ and $\varphi$ for
\[
 W_{u,\varphi}:Q_{p_1}\longrightarrow Q_{p_2}
\]
to be bounded or compact. The parameters $p_1,p_2$ describe the
weights in the derivative norms. Throughout, $p_1>0$ and $0<p_2<1$, and the results are
stated separately for $0<p_1<1$, $p_1=1$, and $p_1>1$.

\subsection{Background and the main question}
Weighted composition operators are often studied by means of measure
embeddings. If $\nu$ is a positive measure and $q>0$, define
\[
 \mu(B)=\int_{\varphi^{-1}(B)}|u(z)|^q\,d\nu(z)
 \qquad(B\subseteq\D\text{ Borel}).
\]
Then
\begin{equation}\label{eq:valuepullbackintro}
 \int_\D|u(z)f(\varphi(z))|^q\,d\nu(z)
 =\int_\D|f(w)|^q\,d\mu(w).
\end{equation}
This identity relates the operator to a Carleson measure problem.
It is used in the Hardy-space results of Contreras and
Hern\'andez-D\'iaz \cite{CH01,CH03} and the Bergman-space results of
\v{C}u\v{c}kovi\'c and Zhao \cite{CZ}. Other descriptions of the
induced measure include Berezin transforms for large Fock spaces
\cite{AT} and sparse domination forms for weighted Bergman spaces \cite{HLSW}.

For spaces defined by derivatives, the product rule gives two terms:
\begin{equation}\label{eq:product}
 (W_{u,\varphi}f)'=\mathcal U f+\mathcal V f,\qquad
 \mathcal U f=u'(f\circ\varphi),\quad
 \mathcal V f=u\varphi'(f'\circ\varphi).
\end{equation}
Pointwise estimates lead to explicit conditions on Bloch and
Bloch-type spaces \cite{OZ,OSZ}. Integral derivative norms, however,
require both terms to be controlled within the same integral; the
same difficulty appears for the Dirichlet space \cite{CGP}. Cheng,
Kumar and Zhou \cite{CKZ} study weighted composition operators on
Dirichlet-type spaces and related $Q_p$ spaces under additional
assumptions on the weight. Boundedness, essential norms and
compactness on BMOA are studied in \cite{Laitila,LL,LLN}.

The identity $Q_1=\mathrm{BMOA}$ relates these questions to the
$Q_p$ scale. Here the operator takes values in $Q_{p_2}$ with
$0<p_2<1$. Since $Q_{p_2}\subseteq\mathrm{BMOA}$, boundedness into
BMOA alone does not give the required estimate, and we keep the
$Q_{p_2}$ norm in both terms of \eqref{eq:product}.

For $p>0$, let
\begin{equation}\label{eq:Nnorm}
 \mathcal N_p=\left\{h\in\Hol(\D):
 \|h\|_{\mathcal N_p}^2
 =\sup_{a\in\D}\int_\D|h(z)|^2
       (1-|\varphi_a(z)|^2)^p\,dA(z)<\infty\right\}.
\end{equation}
These spaces were studied on the disc by Palmberg \cite{Palmberg}
and on the ball by Hu and Khoi \cite{HK}; see also \cite{KM}. With
our norms,
\[
 \|f'\|_{\mathcal N_p}=\|f\|_{Q_p,*},\qquad
 \left\|\int_0^z h(\zeta)\,d\zeta\right\|_{Q_p}
 =\|h\|_{\mathcal N_p}.
\]
We therefore study boundedness and compactness of the two maps in
\eqref{eq:product} from $Q_{p_1}$ to $\mathcal N_{p_2}$.

When $u=1$, only the derivative term remains. Hu and Zhou
\cite{HZcomp} characterize composition operators on $Q_p$ by a
dyadic derivative embedding, and Hu, Luo, Xiao and Zhou \cite{HLXZ}
develop the corresponding theory for $Q_K$ spaces. The
embedding is treated by different capacities in \cite{CJL}.
We use these two embedding theorems as tools; their definitions and
the estimates needed here are given in Section~\ref{sec:prelim}.

The main question is how to pass from an estimate for
$\mathcal U f+\mathcal V f$ to separate estimates for its two terms.
The triangle inequality gives a bound for the sum from separate
bounds, but does not give the converse, because
\[
 |\mathcal U f+\mathcal V f|^2
 =|\mathcal U f|^2+|\mathcal V f|^2
   +2\operatorname{Re}(\mathcal U f\,\overline{\mathcal V f}).
\]
The separation theorem proves
\[
 \begin{aligned}
 W_{u,\varphi}:Q_{p_1}\to Q_{p_2}\text{ bounded}
 &\iff \mathcal U,\mathcal V:Q_{p_1}\to\mathcal N_{p_2}\text{ bounded},\\
 W_{u,\varphi}:Q_{p_1}\to Q_{p_2}\text{ compact}
 &\iff \mathcal U,\mathcal V:Q_{p_1}\to\mathcal N_{p_2}\text{ compact}.
 \end{aligned}
\]
Two analytic functions on each set of a fixed partition give a
uniform estimate for their values and derivatives. Averaging their
finite sums over independent signs gives the estimate for
$\mathcal V$; then $\mathcal U f=(W_{u,\varphi}f)'-\mathcal V f$.
The corresponding estimates near the boundary give compactness.
No further assumption on $\varphi$ is needed.

\subsection{Main results}
We first specify the measures and quantities in the statements. For
$a\in\D$, $0<t<1$, and a Borel set $B\subseteq\D$, put
\begin{align}
 \mu^0_{a,t}(B)&=\int_{\{|z|<t:\,\varphi(z)\in B\}}
 |u'(z)|^2(1-|\varphi_a(z)|^2)^{p_2}\,dA(z),\label{eq:mu0}\\
 \mu^1_a(B)&=\int_{\varphi^{-1}(B)}
 |u(z)\varphi'(z)|^2(1-|\varphi_a(z)|^2)^{p_2}\,dA(z).
 \label{eq:mu1}
\end{align}
The restriction $|z|<t$ makes $\mu^0_{a,t}$ finite; we take the
supremum over $0<t<1$ below.

Section~\ref{sec:prelim} defines the value capacities
$C_{p,\mathcal R}^{(b)}(\nu;E)$ and $D_{p,\mathcal R}^{(b)}(\nu;E)$,
the derivative gauge $\|h\|_{B_p^2(\mathcal D)}$, and
$h_{\nu,2}^p(I)=\nu(Q_I^{\mathrm{up}})/|I|^{p+2}$.
Here $\mathcal D$ is a fixed dyadic grid on $\T$, $\mathcal R$ is the
polar partition used for reconstruction, and $b=b(p_1)$ is fixed as
in \eqref{eq:parameters}. For Borel $E\subseteq\D$, set
\begin{align}
 \mathcal C_{0}^{p_1,p_2}(E)
 &=\sup_{a\in\D,\,0<t<1}
          C_{p_1,\mathcal R}^{(b)}(\mu^0_{a,t};E),\label{eq:Cglobal}\\
 \mathcal C_{1}^{p_1,p_2}(E)
 &=\sup_{a\in\D}
          \|h_{\ind_E\mu^1_a,2}^{p_1}\|_{B_{p_1}^2(\mathcal D)}.
 \nonumber
\end{align}
The subscripts $0$ and $1$ correspond to function values and
derivatives. Both quantities are defined by the symbols and the
explicit coefficient constraints in Section~\ref{sec:prelim}.

Write $E_\rho=\{w\in\D:|w|>\rho\}$, and let $dm$ be normalized
arc measure on $\T$. For $p_1=1$, put
\begin{equation}\label{eq:endpointmax}
 \mathcal M_{p_2}(E)
 =\sup_{a\in\D}\int_\T
       \sup_{\substack{I\in\mathcal D\\\xi\in I}}
       \frac{\mu_a^1(E\cap Q_I^{\mathrm{up}})}{|I|^3}\,dm(\xi).
\end{equation}
Here $Q_I^{\mathrm{up}}$ are the upper half tents defined in
\eqref{eq:tents}. For $p_1>1$, put
\begin{equation}\label{eq:blochJ}
 \mathcal J_{p_2}(E)
 =\sup_{a\in\D}\int_{\varphi^{-1}(E)}
       \frac{|u(z)\varphi'(z)|^2}{(1-|\varphi(z)|^2)^2}
       (1-|\varphi_a(z)|^2)^{p_2}\,dA(z).
\end{equation}
All these quantities are allowed to be infinite.

\begin{theorem}\label{thm:subcritical}
Let $0<p_1<1$ and $0<p_2<1$. Then
$W_{u,\varphi}:Q_{p_1}\to Q_{p_2}$ is bounded if and only if
\[
 u\in Q_{p_2},\qquad
 \mathcal C_{0}^{p_1,p_2}(\D)+\mathcal C_{1}^{p_1,p_2}(\D)<\infty.
\]
Under boundedness, it is compact if and only if
\[
 \lim_{\rho\uparrow1}
 [\mathcal C_{0}^{p_1,p_2}(E_\rho)+\mathcal C_{1}^{p_1,p_2}(E_\rho)]=0.
\]
\end{theorem}

\begin{theorem}\label{thm:endpoint}
Let $p_1=1$ and $0<p_2<1$. Then
$W_{u,\varphi}:\mathrm{BMOA}\to Q_{p_2}$ is bounded if and only if
\[
 u\in Q_{p_2},\qquad
 \mathcal C_{0}^{1,p_2}(\D)+\mathcal M_{p_2}(\D)<\infty.
\]
Under boundedness, it is compact if and only if
\[
 \lim_{\rho\uparrow1}
 [\mathcal C_{0}^{1,p_2}(E_\rho)+\mathcal M_{p_2}(E_\rho)]=0.
\]
In fact, for every Borel set $E\subseteq\D$,
\begin{equation}\label{eq:endpointcapacity}
 \mathcal C_{1}^{1,p_2}(E)=\mathcal M_{p_2}(E).
\end{equation}
\end{theorem}

\begin{theorem}\label{thm:bloch}
Let $p_1>1$ and $0<p_2<1$. Then
$W_{u,\varphi}:Q_{p_1}\to Q_{p_2}$, equivalently
$W_{u,\varphi}:\mathcal B\to Q_{p_2}$, is bounded if and only if
\[
 u\in Q_{p_2},\qquad
 \mathcal C_{0}^{2,p_2}(\D)+\mathcal J_{p_2}(\D)<\infty.
\]
Under boundedness, it is compact if and only if
\[
 \lim_{\rho\uparrow1}
 [\mathcal C_{0}^{2,p_2}(E_\rho)+\mathcal J_{p_2}(E_\rho)]=0.
\]
The value capacity is defined with $p=2$, since $Q_{p_1}=Q_2=\mathcal B$
with equivalent norms for $p_1>1$.
\end{theorem}

The condition $u\in Q_{p_2}$ follows from $W_{u,\varphi}1=u$.
The suprema defining the capacities are taken before the limit
$\rho\uparrow1$. The value capacity $C$ may be replaced by its
equivalent form $D$ in all three statements. For $u=1$, the value
term is zero, and the criteria reduce to the derivative embedding
condition for composition operators. The proofs also give two-sided
norm estimates for $W_{u,\varphi}$; see \eqref{eq:mainnorm}.

The paper is organized as follows. Section~\ref{sec:prelim} collects
the capacity definitions and the auxiliary estimates.
Section~\ref{sec:separation} proves the separation theorem and the
compactness criteria, and then deduces
Theorems~\ref{thm:subcritical}--\ref{thm:bloch}.

We write $A\lesssim B$ if $A\leq CB$ and $A\asymp B$ if both
$A\lesssim B$ and $B\lesssim A$. Unless otherwise stated, $b$ and $M$ are chosen as in
\eqref{eq:parameters} with $p=p_1$.

\section{Preliminaries and lemmas}\label{sec:prelim}
This section has two purposes: to fix the capacity notation used in
the main results and to collect the estimates needed for the
separation argument and the compactness criteria. The two terms in
\eqref{eq:product} lead to the embedding problems
\[
 \int_\D|f|^2\,d\nu
 \quad\text{and}\quad
 \int_\D|f'|^2\,d\nu.
\]
The first uses reconstruction of $f-f(0)$ from $f'$, the second the
distribution of $\nu$ on dyadic tents. Throughout this section
$p>0$ is fixed, and all function-space norms in the lemmas use this
same $p$.

\subsection{The two embedding conditions}
Identify $\T$ with $\mathbb{R}/\mathbb{Z}$. Let
\[
 \mathcal D=\{I_{n,k}:n\geq0,\ 0\leq k<2^n\},\qquad
 I_{n,k}=\{e^{2\pi it}:k2^{-n}\leq t<(k+1)2^{-n}\}.
\]
Arc length is normalized so that $|\T|=1$. For any arc $I\subseteq\T$, define
\begin{equation}\label{eq:tents}
 Q_I=\{re^{i\theta}:e^{i\theta}\in I,\ 1-|I|\leq r<1\},\qquad
 Q_I^{\mathrm{up}}=\{re^{i\theta}:e^{i\theta}\in I,\ 1-|I|\leq r<1-|I|/2\}.
\end{equation}
The sets $Q_I^{\mathrm{up}}$, $I\in\mathcal D$, partition $\D$; see \cite[Section 1]{HZcomp}.

Given a measurable function $F$, set
\begin{equation}\label{eq:TK}
 \|F\|_{T_p}^2=\sup_{I\subseteq\T}
 \int_{Q_I}|F(w)|^2\left(\frac{1-|w|}{|I|}\right)^p dA(w).
\end{equation}
The supremum is over all arcs, including $\T$. Restricting it to $\mathcal D$ gives $\|F\|_{T_p,\mathcal D}$. The standard tent-space characterization, recorded in \cite[(9)--(11)]{CJL}, gives
\begin{equation}\label{eq:boxeq}
 \|F\|_{T_p,\mathcal D}\asymp\|F\|_{T_p},\qquad
 \|f'\|_{T_p}\asymp\|f\|_{Q_p,*}.
\end{equation}

Choose integers $b$ and $M$ once for this $p$, with
\begin{equation}\label{eq:parameters}
 b\geq\max\{b_p,\lceil p\rceil+3\},\qquad M=b+1.
\end{equation}
Choose $b_p$ so that the projection estimates hold for every $b\geq b_p$; see
\cite[Section 2.3]{CJL} and \cite[Lemmas 2.2 and 2.3]{HLXZ}.
Define
\begin{equation}\label{eq:kL}
 dA_b(w)=(b+1)(1-|w|^2)^b dA(w),\qquad
 k_b(\zeta,w)=(1-\zeta\overline w)^{-b-2}.
\end{equation}
The kernel used to reconstruct function values is
\begin{equation}\label{eq:Lclosed}
 L_b(\zeta,w)=\int_0^\zeta k_b(\xi,w)\,d\xi
 =\begin{cases}
 \dfrac{(1-\zeta\overline w)^{-b-1}-1}{(b+1)\overline w},&w\ne0,\\
 \zeta,&w=0.
 \end{cases}
\end{equation}
Now set
\begin{equation}\label{eq:Pdef}
 P_{b,0}F(\zeta)=\int_\D L_b(\zeta,w)F(w)\,dA_b(w).
\end{equation}
By \cite[Lemma 2.12 and Proposition 2.13]{CJL},
\begin{equation}\label{eq:reconstruction}
 \|P_{b,0}F\|_{Q_p}\leq C_{p,b}\|F\|_{T_p},\qquad
 P_{b,0}(f')=f-f(0).
\end{equation}

For the value capacity, use the polar partitions from
\cite[Section 2.2]{CJL}:
\[
 \mathcal R_N=\{R_{N,j,k}:0\leq j,k<2^N\},\qquad N\geq0,
\]
where
\[
 R_{N,j,k}=\{re^{i\theta}:j2^{-N}\leq r^2<(j+1)2^{-N},\
       2\pi k2^{-N}\leq\theta<2\pi(k+1)2^{-N}\}.
\]
The origin is assigned to $R_{N,0,0}$, and the sets are ordered by their indices.
Write $\mathcal R=\{\mathcal R_N\}_{N\geq0}$. For a coefficient vector
$c=(c_R)_{R\in\mathcal R_N}$, define
\begin{equation}\label{eq:Xnorm}
 \|c\|_{X_{p,N}}^2
 =\left\|\sum_{R\in\mathcal R_N}c_R\ind_R\right\|_{T_p}^2
 =\sup_{I\subseteq\T}c^*A_{I,N}^{(p)}c,
\end{equation}
with diagonal matrix
\begin{equation}\label{eq:admissible}
 A_{I,N}^{(p)}=\operatorname{diag}\left(
 \int_{R\cap Q_I}\left(\frac{1-|w|}{|I|}\right)^p dA(w)
 :R\in\mathcal R_N\right).
\end{equation}
Accordingly, $c$ is admissible exactly when $\|c\|_{X_{p,N}}\leq1$.

Let $\nu$ be a finite positive Borel measure and $E\subseteq\D$ a Borel set. Define the holomorphic functions
\begin{equation}\label{eq:restrictedpacket}
 \Psi_{R,E}^{(b)}(\zeta)=\int_{R\cap E}L_b(\zeta,w)\,dA_b(w)
\end{equation}
and the associated Gram matrix
\[
 \Gamma_{\nu,N}^{(b)}(E)(R,S)
 =\int_\D\overline{\Psi_{R,E}^{(b)}(\zeta)}
                      \Psi_{S,E}^{(b)}(\zeta)\,d\nu(\zeta).
\]
In particular,
\[
 c^*\Gamma_{\nu,N}^{(b)}(E)c
 =\int_\D\left|\sum_R c_R\Psi_{R,E}^{(b)}\right|^2d\nu.
\]
We employ the following localized versions of the capacities from
\cite[Definitions 3.1 and 3.3]{CJL}:
\begin{equation}\label{eq:restrictedcapacity}
 C_{p,\mathcal R}^{(b)}(\nu;E)
 =\sup_{N\geq0}\sup_{\|c\|_{X_{p,N}}\leq1}
       c^*\Gamma_{\nu,N}^{(b)}(E)c
\end{equation}
and
\begin{equation}\label{eq:restrictedSDP}
 D_{p,\mathcal R}^{(b)}(\nu;E)
 =\sup_{N\geq0}\inf\left\{\sum_{j=1}^m\lambda_j:
 \Gamma_{\nu,N}^{(b)}(E)\preceq
           \sum_{j=1}^m\lambda_jA_{I_j,N}^{(p)}\right\}.
\end{equation}
Here $m\in\mathbb{N}$, $\lambda_j\geq0$, and each $I_j\subseteq\T$ is an arc.
The ordering $X\preceq Y$ means $v^*(Y-X)v\geq0$ for every complex vector $v$.
The infimum over an empty collection is taken to be $+\infty$.
For $E=\D$, these are the cited capacities, and we omit $E$.
The restriction to a Borel set $E$ is used for compactness.

The value embedding theorem, assembled from \cite[Theorems 3.2, 3.5 and 3.6]{CJL}, asserts
\begin{equation}\label{eq:toolvalue}
 \sup_{\substack{f(0)=0\\\|f\|_{Q_p,*}\leq1}}
       \int_\D|f|^2\,d\nu
 \asymp C_{p,\mathcal R}^{(b)}(\nu)
 \asymp D_{p,\mathcal R}^{(b)}(\nu).
\end{equation}
The form adapted for compactness is
\begin{equation}\label{eq:restrictedembedding}
 \sup_{\|F\|_{T_p}\leq1}
     \int_\D|P_{b,0}(\ind_EF)|^2d\nu
 \asymp C_{p,\mathcal R}^{(b)}(\nu;E)
 \asymp D_{p,\mathcal R}^{(b)}(\nu;E),
\end{equation}
with comparison constants independent of $E$ and $\nu$.
For the first comparison, use the conditional expectations
\[
 E_N^{(b)}F=\sum_{R\in\mathcal R_N}
       \frac{\int_RF\,dA_b}{A_b(R)}\ind_R,
 \qquad A_b(R)=\int_RdA_b.
\]
By \cite[Lemma 2.8 and the proof of Theorem 3.2]{CJL},

\[
 \|E_N^{(b)}F\|_{T_p}\leq C_{p,b}\|F\|_{T_p},\qquad
 \|\ind_E(E_N^{(b)}F-F)\|_{L^1(dA_b)}\longrightarrow0.
\]
For fixed $\zeta$, the kernel $L_b(\zeta,\cdot)$ is bounded.
Fatou's lemma therefore gives
\[
 \int_\D|P_{b,0}(\ind_EF)|^2d\nu
 \leq C_p C_{p,\mathcal R}^{(b)}(\nu;E)\|F\|_{T_p}^2.
\]
The reverse inequality follows by taking $F=\sum_Rc_R\ind_R$.
The second equivalence is a consequence of the proof of
\cite[Theorem 3.5]{CJL}: its matrix argument applies to every positive
semidefinite matrix with the same diagonal constraints.
Apply it to $\Gamma_{\nu,N}^{(b)}(E)$ at each level $N$, and then
take the supremum over $N$. The constants do not depend on $E$.

For the derivative embedding, we use the capacity gauge from
\cite[(1.2)]{HLXZ}. Given $h:\mathcal D\to[0,\infty]$, define
\begin{equation}\label{eq:treegauge}
 \|h\|_{B_p^2(\mathcal D)}=\sup_d\sum_{I\in\mathcal D}d_Ih(I),
\end{equation}
where the admissible weights satisfy
\begin{equation}\label{eq:treeadmissible}
 d_I\geq0,\qquad
 \sum_{I\subseteq J}d_I\leq |J|^p\quad(J\in\mathcal D).
\end{equation}
All sums are over nonnegative terms, with the convention $0\cdot\infty=0$.
No square root is taken in this definition. For a positive Borel
measure $\nu$, set
\begin{equation}\label{eq:hmeasure}
 h_{\nu,2}^p(I)=\frac{\nu(Q_I^{\mathrm{up}})}{|I|^{p+2}}.
\end{equation}
The derivative embedding theorem gives
\begin{equation}\label{eq:localtrace}
 \sup_{\|f\|_{Q_p,*}\leq1}\int_\D|f'|^2\,d\nu
 \asymp\|h_{\nu,2}^p\|_{B_p^2(\mathcal D)}.
\end{equation}
For absolutely continuous measures this is \cite[Theorem 1.1]{HLXZ}, with exponent two.
The proof after Lemma~\ref{lem:atoms} gives the same result for positive Borel measures. The constants are unchanged under restriction of $\nu$.

The two restrictions have different meanings:
\[
 C_{p,\mathcal R}^{(b)}(\nu;E):\quad F\mapsto P_{b,0}(\ind_EF),
 \qquad
 h_{\ind_E\nu,2}^p:\quad \nu\mapsto\ind_E\nu.
\]

\subsection{Local estimates and finite sums}
We introduce the positive integral operator
\[
 B_b^+F(\zeta)=\int_\D
       \frac{|F(w)|}{|1-\zeta\overline w|^{b+2}}\,dA_b(w).
\]
For $b$ as in \eqref{eq:parameters}, the absolute-kernel estimate in \cite[Lemma 2.2 and the proof of
Lemma 2.3]{HLXZ} gives
\begin{equation}\label{eq:positiveprojection}
 \|B_b^+F\|_{T_p}\leq C_{p,b}\|F\|_{T_p}.
\end{equation}
Since $b>p$ and $A_b(\D)=1$, we further obtain
\begin{equation}\label{eq:L1L2}
 \|F\|_{L^1(dA_b)}\leq\|F\|_{L^2(dA_b)}
 \leq C_{p,b}\|F\|_{T_p},
\end{equation}
along with the standard growth estimates recorded in \cite{WZ,Zhu}
\begin{equation}\label{eq:Bloch}
 (1-|z|^2)|f'(z)|\leq C_p\|f\|_{Q_p,*},\qquad
 |f(z)|\leq C_p\log\frac e{1-|z|^2}\|f\|_{Q_p}.
\end{equation}
The norm of the point-evaluation functional satisfies
\begin{equation}\label{eq:eval}
 \|\operatorname{ev}_z\|_{Q_p^*}\asymp_p\log\frac e{1-|z|^2}.
\end{equation}

For the lower bound, use the logarithm
$g(w)=\log(1/(1-w))$, normalized so that $g(0)=0$. Given an arc $I$ with length $h=|I|$,
\[
 \int_{Q_I}\frac{(1-|w|)^p}{|1-w|^2}\,dA(w)
 \leq C\int_{1-h}^1(1-r)^p
       \left(\int_{-\pi}^{\pi}\frac{d\theta}{|1-re^{i\theta}|^2}\right)r\,dr
 \leq C\int_{1-h}^1(1-r)^{p-1}\,dr
 =\frac Cp h^p.
\]
The angular integral equals $2\pi/(1-r^2)$. Applying \eqref{eq:boxeq}, we obtain $\|g\|_{Q_p,*}\leq C_p$.
For any nonzero $z\in\D$, define
\[
 g_z(w)=g\!\left(\frac{\overline z}{|z|}\varphi_z(w)\right)-g(|z|).
\]
M\"obius invariance then implies
\[
 g_z(0)=0,\qquad \|g_z\|_{Q_p}\leq C_p,\qquad
 |g_z(z)|=\log\frac1{1-|z|}.
\]
Combining this family with the constant function $1$, we obtain
\[
 \|\operatorname{ev}_z\|_{Q_p^*}
 \geq\max\left\{1,C_p^{-1}\log\frac1{1-|z|}\right\}
 \gtrsim_p\log\frac e{1-|z|^2}.
\]
The case $z=0$ is immediate, and \eqref{eq:Bloch} gives the upper bound.

For each integer $j\geq0$ and $0\leq l<2^{j+3}$, define
\begin{align}
 V_{j,l}&=\left\{re^{i\theta}:1-2^{-j}\leq r<1-2^{-j-1},
              \ \frac{2\pi l}{2^{j+3}}\leq\theta<\frac{2\pi(l+1)}{2^{j+3}}\right\},
              \label{eq:whitney}\\
 a_{j,l}&=\left(1-\frac{3}{4}\cdot 2^{-j}\right)\exp\!\left(\frac{2\pi i(l+1/2)}{2^{j+3}}\right).
              \label{eq:centers}
\end{align}
Assign the origin to $V_{0,0}$ and use half-open boundaries, so that the sets form a Borel partition of $\D$. Write $\lambda=(j,l)$, $V_\lambda=V_{j,l}$, and $a_\lambda=a_{j,l}$.

The following bounds will be used for this partition.

\begin{lemma}\label{lem:whitneygeometry}
Let $\lambda=(j,l)$. For every $w\in V_{j,l}$,
\begin{gather}
 \frac{15}{16}2^{-j}\leq1-|a_{j,l}|^2\leq\frac32\,2^{-j},\qquad
 \frac{2^{-2j}}{32}\leq A(V_{j,l})\leq\frac{2^{-2j}}8,
 \label{eq:whitneygeometry}\\
 \frac12(1-|a_\lambda|^2)\leq|1-\overline a_\lambda w|
                       \leq4(1-|a_\lambda|^2),\qquad |a_\lambda|\geq\frac14.
 \label{eq:localgeometry}
\end{gather}
For arbitrary $\zeta\in\D$,
\begin{equation}\label{eq:denomcompare}
 \frac13|1-\overline a_\lambda\zeta|
 \leq|1-\overline w\zeta|\leq2|1-\overline a_\lambda\zeta|.
\end{equation}
If we set
\[
 V_\lambda^*=\{z:\operatorname{dist}(z,V_\lambda)
                         <(1-|a_\lambda|^2)/16\},
\]
then $V_\lambda^*\subset\D$, the overlap multiplicity of the enlarged sets is bounded by an absolute constant, and $1-|z|\asymp1-|a_\lambda|^2$ uniformly on $V_\lambda^*$.
\end{lemma}
\begin{proof}
Take $w=re^{i\theta}\in V_{j,l}$, with $\theta$ lying in the corresponding angular sector. Direct calculation gives
\begin{align*}
 1-|a_{j,l}|^2&=\frac34\,2^{-j}\left(2-\frac34\,2^{-j}\right),\\
 A(V_{j,l})&=2^{-j-3}
       \left[(1-2^{-j-1})^2-(1-2^{-j})^2\right]
       =\frac{2^{-2j}}8\left(1-\frac34\,2^{-j}\right),\\
 |w-a_{j,l}|&\leq
 \bigl||w|-|a_{j,l}|\bigr|+|\theta-\arg a_{j,l}|
 \leq\left(\frac14+\frac\pi8\right)2^{-j}.
\end{align*}
Angular differences are measured inside the sector defining $V_{j,l}$.
The estimates \eqref{eq:whitneygeometry} follow immediately, and
\begin{align*}
 |1-\overline a_{j,l}w|&\geq1-|a_{j,l}|=\frac34\,2^{-j}
                   \geq\frac12(1-|a_{j,l}|^2),\\
 |1-\overline a_{j,l}w|&\leq1-|a_{j,l}|^2+|a_{j,l}-w|
                   <3\cdot2^{-j}<4(1-|a_{j,l}|^2).
\end{align*}
Moreover, $|a_{j,l}|\geq1/4$. For any $\zeta\in\D$,
\begin{gather*}
 |1-\overline w\zeta|\geq2^{-j-1},\qquad
 |1-\overline a_{j,l}\zeta|\geq\frac34\,2^{-j},\\
 \bigl||1-\overline w\zeta|-|1-\overline a_{j,l}\zeta|\bigr|
 \leq\left(\frac14+\frac\pi8\right)2^{-j}.
\end{gather*}
We then obtain
\begin{align*}
 |1-\overline w\zeta|
 &\leq\left[1+\frac43\left(\frac14+\frac\pi8\right)\right]
                      |1-\overline a_{j,l}\zeta|
 <2|1-\overline a_{j,l}\zeta|,\\
 |1-\overline a_{j,l}\zeta|
 &\leq\left[1+2\left(\frac14+\frac\pi8\right)\right]
                      |1-\overline w\zeta|
 <3|1-\overline w\zeta|.
\end{align*}
For $z\in V_{j,l}^*$, the dilation radius is bounded above by
$3\cdot2^{-j}/32$, so
\[
 \frac{13}{32}\,2^{-j}<1-|z|<\frac{35}{32}\,2^{-j}.
\]
For a fixed point $z$, this restricts $j$ to a finite set of values.
At each admissible scale, the angular extension has length at most $C2^{-j}$,
while the original sectors have angular width $2\pi2^{-j-3}$. Consequently
\[
 \#\{j:\exists l,\ z\in V_{j,l}^*\}\leq C,\qquad
 \sup_j\#\{l:z\in V_{j,l}^*\}\leq C,
 \qquad \sum_{j,l}\ind_{V_{j,l}^*}(z)\leq C.
\]
The level $j=0$ has only eight sectors. This proves all assertions.
\end{proof}

Given a coefficient sequence $c=(c_\lambda)$, define
\begin{equation}\label{eq:seqnorm}
 F_c(w)=\sum_\lambda\frac{c_\lambda}{(1-|a_\lambda|^2)}\ind_{V_\lambda}(w),
 \qquad \|c\|_{\mathfrak t_p}=\|F_c\|_{T_p}.
\end{equation}
This norm depends solely on the moduli $|c_\lambda|$ and does not increase upon removing coordinates from the sequence. We also have
\begin{equation}\label{eq:coeffsup}
 \sup_\lambda|c_\lambda|\leq C_p\|c\|_{\mathfrak t_p}.
\end{equation}
To justify \eqref{eq:coeffsup}, take an arc containing the angular projection of $V_\lambda$, with length $\min\{8(1-|a_\lambda|^2),1\}$. Its tent contains $V_\lambda$, and on this set
\[
 \left(\frac{1-|w|}{|I|}\right)^p\geq 24^{-p}.
\]
An application of Lemma~\ref{lem:whitneygeometry} gives
\[
 \|c\|_{\mathfrak t_p}^2
 \geq\frac{|c_\lambda|^2A(V_\lambda)}{(1-|a_\lambda|^2)^2}24^{-p}
 \geq c_p|c_\lambda|^2.
\]
For later reference, pairwise disjointness of the sets gives the exact representation
\begin{equation}\label{eq:seqexplicit}
 \|c\|_{\mathfrak t_p}^2
 =\sup_I\sum_\lambda\frac{|c_\lambda|^2}{(1-|a_\lambda|^2)^2}
       \int_{V_\lambda\cap Q_I}\left(\frac{1-|w|}{|I|}\right)^p dA(w).
\end{equation}
It follows that for any index set $\Lambda$ and any choice of unimodular scalars $\theta_\lambda$,
\[
 \|(c_\lambda\ind_\Lambda(\lambda))\|_{\mathfrak t_p}
 \leq\|c\|_{\mathfrak t_p},\qquad
 \|(\theta_\lambda c_\lambda)\|_{\mathfrak t_p}=\|c\|_{\mathfrak t_p}.
\]

\begin{lemma}\label{lem:sampling}
For $f\in Q_{p}$, set
\begin{equation}\label{eq:samplecoeff}
 c_\lambda(f)=(1-|a_\lambda|^2)\sup_{w\in V_\lambda}|f'(w)|.
\end{equation}
Then $\|c(f)\|_{\mathfrak t_p}\leq C_p\|f\|_{Q_{p},*}$.
\end{lemma}
\begin{proof}
For $w\in V_\lambda$, the disk
$D(w,(1-|a_\lambda|^2)/32)$ is contained in $V_\lambda^*$.
The submean property of subharmonic functions gives
\[
 \sup_{V_\lambda}|f'|^2\leq
 \frac{C}{(1-|a_\lambda|^2)^2}\int_{V_\lambda^*}|f'(\zeta)|^2\,dA(\zeta).
\]
Fix an arc $I$ and write $h=|I|$. Let $J$ denote the concentric arc with length $\min\{16h,1\}$, and take $J=\T$ whenever $16h\geq1$.
From the definitions of $V_\lambda$ and $V_\lambda^*$,
\begin{gather*}
 V_\lambda\cap Q_I\ne\varnothing
 \ \Longrightarrow\quad
 1-|a_\lambda|^2\leq3h,\quad V_\lambda^*\subseteq Q_J,\\
 \left(\frac{1-|w|}{h}\right)^p
 \leq C_p\left(\frac{1-|\zeta|}{|J|}\right)^p
 \quad(w\in V_\lambda\cap Q_I,\ \zeta\in V_\lambda^*).
\end{gather*}
Using the bound $A(V_\lambda)\lesssim(1-|a_\lambda|^2)^2$, we compute
\begin{align*}
 &\int_{Q_I}|F_{c(f)}(w)|^2
                 \left(\frac{1-|w|}{h}\right)^p dA(w)\\
 &=\sum_{V_\lambda\cap Q_I\ne\varnothing}
       \sup_{V_\lambda}|f'|^2
       \int_{V_\lambda\cap Q_I}\left(\frac{1-|w|}{h}\right)^p dA(w)\\
 &\leq C_p\sum_{V_\lambda\cap Q_I\ne\varnothing}
       \frac{A(V_\lambda\cap Q_I)}{(1-|a_\lambda|^2)^2}
       \int_{V_\lambda^*}|f'(\zeta)|^2
                    \left(\frac{1-|\zeta|}{|J|}\right)^p dA(\zeta)\\
 &\leq C_p\int_{Q_J}|f'(\zeta)|^2
                    \left(\frac{1-|\zeta|}{|J|}\right)^p
                    \sum_\lambda\ind_{V_\lambda^*}(\zeta)\,dA(\zeta)\\
 &\leq C_p\|f'\|_{T_p}^2\leq C_p\|f\|_{Q_{p},*}^2.
\end{align*}
Taking the supremum over $I$ proves the estimate.
\end{proof}

\begin{lemma}\label{lem:atoms}
Let $M$ be as in \eqref{eq:parameters}. For $q=M,M+1$, define
\begin{equation}\label{eq:atoms}
 h_{\lambda,q}(\zeta)=\left(\frac{(1-|a_\lambda|^2)}{1-\overline a_\lambda\zeta}\right)^q,
 \qquad S_qc=\sum_\lambda c_\lambda h_{\lambda,q}
\end{equation}
on the space of finitely supported sequences $c$. Then
\begin{equation}\label{eq:atomsbound}
 \|S_qc\|_{Q_{p}}\leq C_{p,M}\|c\|_{\mathfrak t_p}.
\end{equation}
If the support of $c$ is confined to indices with $j\geq J$, then for every $R<1$ and integer $k\geq0$,
\begin{equation}\label{eq:atomstail}
 \sup_{|\zeta|\leq R}|(S_qc)^{(k)}(\zeta)|
 \leq C_{p,M,R,k}\,2^{-J(q-1)}\|c\|_{\mathfrak t_p}.
\end{equation}
\end{lemma}
\begin{proof}
For $w\in V_\lambda$, Lemma~\ref{lem:whitneygeometry} guarantees
\[
 1-|w|^2\asymp1-|a_\lambda|^2,\qquad
 A(V_\lambda)\asymp(1-|a_\lambda|^2)^2.
\]
Accordingly, the positive integral over a single set satisfies
\begin{align*}
 \int_{V_\lambda}
     \frac{|F_c(w)|(1-|w|^2)^{q-1}}{|1-\overline w\zeta|^{q+1}}\,dA(w) \asymp
     \frac{|c_\lambda|(1-|a_\lambda|^2)^q}
          {|1-\overline a_\lambda\zeta|^{q+1}}.
\end{align*}
Since
\[
 h'_{\lambda,q}(\zeta)=
 \frac{q\overline a_\lambda(1-|a_\lambda|^2)^q}
      {(1-\overline a_\lambda\zeta)^{q+1}},
\]
summing these estimates gives
\begin{align*}
 |(S_qc)'(\zeta)|
 &\leq q\sum_\lambda
      \frac{|c_\lambda|(1-|a_\lambda|^2)^q}
           {|1-\overline a_\lambda\zeta|^{q+1}}\\
 &\leq C_q\int_\D
       \frac{|F_c(w)|(1-|w|^2)^{q-1}}{|1-\overline w\zeta|^{q+1}}\,dA(w)
 \leq C_qB_{q-1}^+F_c(\zeta).
\end{align*}
Applying \eqref{eq:positiveprojection} with exponent $q-1$,
\[
 \|S_qc\|_{Q_{p},*}
 \lesssim\|(S_qc)'\|_{T_p}
 \lesssim\|B_{q-1}^+F_c\|_{T_p}
 \lesssim\|F_c\|_{T_p}=\|c\|_{\mathfrak t_p}.
\]
At the origin, \eqref{eq:coeffsup} gives
\begin{align*}
 |S_qc(0)|
 &\leq\sup_\lambda|c_\lambda|
       \sum_{j\geq0}\sum_{l=0}^{2^{j+3}-1}(1-|a_{j,l}|^2)^q\\
 &\leq C_q\sup_\lambda|c_\lambda|\sum_{j\geq0}2^{j+3}2^{-jq}
 \leq C_{p,q}\|c\|_{\mathfrak t_p}.
\end{align*}
Because $q>1$, the series converges absolutely. The same argument shows that for any unimodular $\theta_\lambda$,
\[
 \left\|\sum_\lambda\theta_\lambda c_\lambda h_{\lambda,q}\right\|_{Q_{p}}
 \leq C_{p,M}\|c\|_{\mathfrak t_p}.
\]
The constant does not depend on the finite index set or on the choice of phases $\theta_\lambda$.

For $|\zeta|\leq R$ and $k\geq0$,
\[
 |h_{\lambda,q}^{(k)}(\zeta)|
 \leq C_{q,k}(1-R)^{-q-k}(1-|a_\lambda|^2)^q.
\]
Hence, when $\supp c\subseteq\{j\geq J\}$,
\begin{align*}
 \sup_{|\zeta|\leq R}|(S_qc)^{(k)}(\zeta)|
 &\leq C_{q,R,k}\sup_\lambda|c_\lambda|
                  \sum_{j\geq J}2^{j+3}2^{-jq}\\
 &\leq C_{p,q,R,k}\|c\|_{\mathfrak t_p}
            \frac{2^{-J(q-1)}}{1-2^{1-q}}.
\end{align*}
This proves \eqref{eq:atomstail}.
\end{proof}

The coefficient norm introduced in \eqref{eq:seqnorm} gives another formula for the derivative gauge:
\begin{gather}
 \mathfrak L_p(\nu)=\sup_{\substack{c\text{ finite}\\\|c\|_{\mathfrak t_p}\leq1}}
 \sum_\lambda\frac{\nu(V_\lambda)}{(1-|a_\lambda|^2)^2}|c_\lambda|^2,
 \label{eq:coefficientgauge}\\
 \mathfrak L_p(\nu)\asymp\|h_{\nu,2}^p\|_{B_p^2(\mathcal D)}.\label{eq:treecomparison}
\end{gather}
Each upper half $Q_I^{\mathrm{up}}$ contains exactly eight sets $V_\lambda$.
Whenever $V_\lambda\subset Q_I^{\mathrm{up}}$ and $I\subseteq J$, we have
$1-|a_\lambda|^2\asymp|I|$, $1-|w|\asymp|I|$ uniformly on $V_\lambda$, and
$A(V_\lambda)\asymp|I|^2$. It follows that
\[
 \int_{V_\lambda}\frac{1}{(1-|a_\lambda|^2)^2}
       \left(\frac{1-|w|}{|J|}\right)^p dA(w)
 \asymp_p\frac{|I|^2}{|I|^2}
       \left(\frac{|I|}{|J|}\right)^p
 =\left(\frac{|I|}{|J|}\right)^p.
\]
For a dyadic interval $J$, its tent decomposes disjointly into the upper halves indexed by all $I\subseteq J$. Summing this estimate and invoking
\eqref{eq:boxeq} gives
\begin{equation}\label{eq:general-lattice}
 \|c\|_{\mathfrak t_p}^2\asymp_p
 \sup_{J\in\mathcal D}\sum_{I\subseteq J}(|I|/|J|)^p
             \sum_{V_\lambda\subset Q_I^{\mathrm{up}}}|c_\lambda|^2.
\end{equation}
Given a finite coefficient family, set
$d_I=|I|^p\max_{V_\lambda\subset Q_I^{\mathrm{up}}}|c_\lambda|^2$.
Then
\begin{gather*}
 \sup_J\sum_{I\subseteq J}|J|^{-p}d_I
 \lesssim_p\|c\|_{\mathfrak t_p}^2,\\
 \sum_\lambda\frac{\nu(V_\lambda)|c_\lambda|^2}{(1-|a_\lambda|^2)^2}
 \lesssim\sum_I d_I\frac{\nu(Q_I^{\mathrm{up}})}{|I|^{p+2}}.
\end{gather*}
Conversely, for an admissible finite sequence $d$, assign
$c_\lambda=(d_I/|I|^p)^{1/2}$ to each of the eight cells inside $Q_I^{\mathrm{up}}$.
One obtains
\begin{gather*}
 \|c\|_{\mathfrak t_p}^2
 \lesssim_p 8\sup_J\sum_{I\subseteq J}
       |J|^{-p}d_I\leq C_p',\\
 \sum_\lambda\frac{\nu(V_\lambda)|c_\lambda|^2}{(1-|a_\lambda|^2)^2}
 \asymp\sum_I d_I\frac{\nu(Q_I^{\mathrm{up}})}{|I|^{p+2}}.
\end{gather*}
After a fixed rescaling of coefficients, these inequalities prove
\eqref{eq:treecomparison}. The constants are independent of $\nu$ and of any measurable restriction of $\nu$.

Also, Lemmas~\ref{lem:sampling} and \ref{lem:atoms} together imply \eqref{eq:localtrace} for positive Borel measures. Define
\[
 E_p(\nu)=\sup_{\|f\|_{Q_p,*}\leq1}\int_\D|f'|^2d\nu.
\]
A sampling argument combined with finite truncation and monotone convergence
gives $E_p(\nu)\leq C_p\mathfrak L_p(\nu)$. For the reverse inequality, note that
$|h'_{\lambda,M}(w)|\geq c_M/(1-|a_\lambda|^2)$ on $V_\lambda$.
Take a finitely supported sequence $c$ and independent sign variables $\varepsilon_\lambda=\pm1$. Then
\begin{align*}
 \sum_\lambda\frac{|c_\lambda|^2\nu(V_\lambda)}{(1-|a_\lambda|^2)^2}
 &\leq C_M\sum_\lambda |c_\lambda|^2
                 \int_\D |h'_{\lambda,M}|^2d\nu\\
 &=C_M\mathbb E\int_\D|(S_M(\varepsilon c))'|^2d\nu\\
 &\leq C_{p,M}E_p(\nu)\|c\|_{\mathfrak t_p}^2.
\end{align*}
Hence $\mathfrak L_p(\nu)\leq C_pE_p(\nu)$.
Together with \eqref{eq:treecomparison}, this proves the equivalence. The expectation involves only finite sums; truncation handles infinite integrals. The constants remain unchanged under measurable restrictions of $\nu$.

\begin{lemma}\label{lem:jetmatrix}
For $w\in V_\lambda$ and arbitrary $x,y\in\C$,
\begin{equation}\label{eq:jetseparate}
 |x|^2+(1-|a_\lambda|^2)^{-2}|y|^2
 \leq C_M\sum_{q=M}^{M+1}
          |x h_{\lambda,q}(w)+y h'_{\lambda,q}(w)|^2.
\end{equation}
\end{lemma}
\begin{proof}
For $w\in V_\lambda$, Lemma~\ref{lem:whitneygeometry} gives
\[
 4^{-q}\leq|h_{\lambda,q}(w)|\leq2^q,
 \qquad q=M,M+1,
\]
and
\[
 \frac1{16}\leq
 \left|\frac{(1-|a_\lambda|^2)\overline a_\lambda}
 {1-\overline a_\lambda w}\right|\leq2.
\]
Compute the logarithmic derivative:
\[
 \frac{h'_{\lambda,q}(w)}{h_{\lambda,q}(w)}
 =\frac{q\overline a_\lambda}{1-\overline a_\lambda w}.
\]
We then derive the identities
\begin{align*}
 \frac{\overline a_\lambda y}{1-\overline a_\lambda w}
 &=\frac{xh_{\lambda,M+1}(w)+yh'_{\lambda,M+1}(w)}{h_{\lambda,M+1}(w)}
   -\frac{xh_{\lambda,M}(w)+yh'_{\lambda,M}(w)}{h_{\lambda,M}(w)},\\
 x&=(M+1)\frac{xh_{\lambda,M}(w)+yh'_{\lambda,M}(w)}{h_{\lambda,M}(w)}
   -M\frac{xh_{\lambda,M+1}(w)+yh'_{\lambda,M+1}(w)}{h_{\lambda,M+1}(w)}.
\end{align*}
Using the first identity and the lower bound $1/16$, we obtain
\[
 \frac{|y|}{1-|a_\lambda|^2}
 \leq16\sum_{q=M}^{M+1}4^q
      |xh_{\lambda,q}(w)+yh'_{\lambda,q}(w)|.
\]
The second identity provides
\[
 |x|\leq(M+1)\sum_{q=M}^{M+1}4^q
      |xh_{\lambda,q}(w)+yh'_{\lambda,q}(w)|.
\]
Square both estimates and apply the convexity bound $(s+t)^2\leq2s^2+2t^2$.
Adding the resulting inequalities gives \eqref{eq:jetseparate}.
\end{proof}

\subsection{Compactness and reconstruction}
Let
\[
 E_\rho=\{w\in\D:|w|>\rho\},\qquad
 P_\rho F=P_{b,0}\bigl(\ind_{E_\rho}F\bigr),\qquad 0<\rho<1.
\]
We introduce the positive quantity
\begin{equation}\label{eq:eta}
 \eta_b(\rho)=
 \left(\int_{|w|>\rho}\frac{(1-|w|^2)^{2b}}{(1-|w|)^p}\,dA(w)\right)^{1/2}.
\end{equation}
A direct integration gives the estimate
\begin{equation}\label{eq:etarate}
 \eta_b(\rho)\leq C_{p,b}(1-\rho)^{b+(1-p)/2}\longrightarrow0.
\end{equation}
The exponent appearing above is positive under the standing choice of $b$.

\begin{lemma}\label{lem:taillocal}
There exists a constant such that, uniformly over all $F$ satisfying $\|F\|_{T_p}\leq1$,
\begin{equation}\label{eq:tailbounded}
 \|P_\rho F\|_{Q_{p}}\leq C_{p,b},
\end{equation}
and, for each fixed $0<R<1$,
\begin{equation}\label{eq:taillocal}
 \sup_{|\zeta|\leq R}
 \bigl(|P_\rho F(\zeta)|+|(P_\rho F)'(\zeta)|\bigr)
 \leq C_{p,b,R}\eta_b(\rho).
\end{equation}
\end{lemma}
\begin{proof}
From the contractive bound $\|\ind_{E_\rho}F\|_{T_p}\leq\|F\|_{T_p}$ together with identity \eqref{eq:reconstruction},
\[
 \|P_\rho F\|_{Q_{p}}
 \leq C_{p,b}\|\ind_{E_\rho}F\|_{T_p}
 \leq C_{p,b}\|F\|_{T_p},\qquad P_\rho F(0)=0.
\]
For any $R<1$, the kernel admits the uniform bound
\[
 \sup_{\substack{|\zeta|\leq R\\w\in\D}}
       \bigl(|L_b(\zeta,w)|+|k_b(\zeta,w)|\bigr)\leq C_{b,R}.
\]
An application of the Cauchy–Schwarz inequality, combined with the definition of the root tent, gives
\begin{align*}
 &\sup_{|\zeta|\leq R}\bigl(|P_\rho F(\zeta)|+|(P_\rho F)'(\zeta)|\bigr)\\
 &\leq C_{b,R}\int_{E_\rho}|F(w)|\,dA_b(w)\\
 &\leq C_{b,R}
  \left(\int_{E_\rho}|F(w)|^2(1-|w|)^p\,dA(w)\right)^{1/2}
  \left(\int_{E_\rho}\frac{(1-|w|^2)^{2b}}{(1-|w|)^p}\,dA(w)\right)^{1/2}\\
 &\leq C_{p,b,R}\|F\|_{T_p}\eta_b(\rho).
\end{align*}
We next examine the asymptotic behaviour of $\eta_b$. Changing variables to radial coordinates,
\begin{align*}
 \eta_b(\rho)^2
 &=2\int_\rho^1(1-r)^{2b-p}(1+r)^{2b}r\,dr\\
 &\leq\frac{2^{2b+1}}{2b-p+1}(1-\rho)^{2b-p+1}
 \longrightarrow0,
\end{align*}
where convergence follows from the hypothesis $2b-p+1>0$.
\end{proof}

We shall also need the part of reconstruction over $|w|\leq\rho$.

\begin{lemma}\label{lem:interiorcompact}
For $0<\rho<1$, the mapping
\begin{equation}\label{eq:Rrho}
 R_\rho F=P_{b,0}\bigl(\ind_{\{|w|\leq\rho\}}F\bigr)
\end{equation}
is compact as an operator from $T_p$ into $Q_{p}$. Consequently,
\begin{equation}\label{eq:Srho}
 S_\rho f=f(0)+R_\rho(f')
\end{equation}
defines a compact operator on $Q_{p}$, whose operator norm remains bounded uniformly in $\rho$.
\end{lemma}
\begin{proof}
The reproducing kernel admits the series expansion
\[
 k_b(\zeta,w)=\sum_{j=0}^\infty a_{b,j}\zeta^j\overline w^j,
 \qquad a_{b,j}=\binom{b+j+1}{j}\leq C_b(j+1)^{b+1}.
\]
Introduce the finite-rank approximation
\[
 R_{\rho,N}F(\zeta)=\sum_{j=0}^N\frac{a_{b,j}}{j+1}\zeta^{j+1}
             \int_{|w|\leq\rho}F(w)\overline w^j\,dA_b(w).
\]
For $\|F\|_{T_p}\leq1$, inequality \eqref{eq:L1L2} gives
\[
 \left|a_{b,j}\int_{|w|\leq\rho}F(w)\overline w^j\,dA_b(w)\right|
 \leq C_{p,b}(j+1)^{b+1}\rho^j.
\]
If $g$ vanishes at the origin, then
\begin{equation}\label{eq:supderivative}
 \|g\|_{Q_{p}}^2
 =\sup_a\int_\D|g'|^2(1-|\varphi_a|^2)^{p}\,dA
 \leq\|g'\|_\infty^2,
\end{equation}
which implies
\begin{align*}
 \|R_\rho-R_{\rho,N}\|_{T_p\to Q_{p}}
 &\leq\sup_{\|F\|_{T_p}\leq1}
                   \|(R_\rho F-R_{\rho,N}F)'\|_\infty\\
 &\leq C_{p,b}\sum_{j>N}(j+1)^{b+1}\rho^j\longrightarrow0.
\end{align*}
Thus $R_\rho$ is an operator-norm limit of finite-rank operators and is compact. The decomposition
\[
 S_\rho f=f(0)+R_\rho(f'),\qquad
 \|S_\rho f\|_{Q_{p}}
 \leq |f(0)|+C_p\|f'\|_{T_p}\leq C_p\|f\|_{Q_{p}}
\]
proves compactness and uniform boundedness of $S_\rho$.
\end{proof}

The reconstruction formula gives the exact identity
\begin{equation}\label{eq:decompidentity}
 f-S_\rho f=P_\rho(f').
\end{equation}

\subsection{The endpoint maximal identity}
When $p=1$, the coefficient condition \eqref{eq:treeadmissible} reduces to
$\sum_{I\subseteq J}d_I\leq |J|$. We give a direct proof of its maximal integral representation. Set
\[
 \mathcal D_N=\{I\in\mathcal D:|I|\geq2^{-N}\},\qquad N\geq0.
\]

\begin{lemma}\label{lem:endpointpacking}
For every function $h:\mathcal D\to[0,\infty]$,
\begin{equation}\label{eq:packingmax}
 \sup_{\substack{d_I\geq0\\
        \sum_{I\subseteq J}d_I\leq|J|\ (J\in\mathcal D)}}
       \sum_{I\in\mathcal D}d_Ih(I)
 =\int_\T\sup_{\substack{I\in\mathcal D\\\xi\in I}}h(I)\,dm(\xi).
\end{equation}
Both sides may be infinite.
\end{lemma}

\begin{proof}
Fix an integer $N\geq0$ and a real number $R>0$, and define
\[
 h_R(I)=\min\{h(I),R\},\qquad
 M_{N,R}(\xi)=\max_{\substack{I\in\mathcal D_N\\\xi\in I}}h_R(I).
\]
The maximum is attained because the family is finite and contains the root interval $\T$.
Let
\[
 P_{N,R}=\sup\left\{
    \sum_{I\in\mathcal D_N}d_Ih_R(I):
    d_I\geq0,\quad
    \sum_{\substack{I\in\mathcal D_N\\I\subseteq J}}d_I\leq |J|
    \quad(J\in\mathcal D_N)\right\}.
\]
We first verify
\begin{equation}\label{eq:endpointfinite}
 P_{N,R}=\int_\T M_{N,R}\,dm.
\end{equation}

Take any admissible coefficients $(d_I)_{I\in\mathcal D_N}$. There exist pairwise disjoint measurable subsets $E_I\subseteq I$ satisfying
\begin{equation}\label{eq:endpointallocation}
 m(E_I)=d_I\qquad(I\in\mathcal D_N).
\end{equation}
These sets may be constructed inductively. Begin with arcs of generation $N$, selecting $E_I\subseteq I$ with measure $d_I\leq |I|$. Suppose the sets for all arcs of generation strictly larger than $n$ have already been chosen. For an arc $J$ at generation $n$, define
\[
 F_J=J\setminus
       \bigcup_{\substack{I\in\mathcal D_N\\I\subsetneq J}}E_I.
\]
The previously selected sets are pairwise disjoint, so
\[
 m(F_J)=|J|-\sum_{\substack{I\in\mathcal D_N\\I\subsetneq J}}d_I
       \geq d_J.
\]
Since arc length measure is non-atomic, we can choose a measurable subset $E_J\subseteq F_J$ with measure exactly $d_J$. Arcs of a fixed generation are disjoint, and sets chosen outside $J$ do not intersect $J$. This completes the induction and gives \eqref{eq:endpointallocation}.

For $\xi\in E_I\subseteq I$, we have $h_R(I)\leq M_{N,R}(\xi)$. It follows that
\begin{align*}
 \sum_{I\in\mathcal D_N}d_Ih_R(I)
 &=\sum_{I\in\mathcal D_N}\int_{E_I}h_R(I)\,dm\\
 &\leq\sum_{I\in\mathcal D_N}\int_{E_I}M_{N,R}\,dm
 \leq\int_\T M_{N,R}\,dm.
\end{align*}
Taking the supremum gives the upper bound in \eqref{eq:endpointfinite}.

For the reverse inequality, enumerate
$\mathcal D_N=\{I_1,\ldots,I_L\}$. Put
\[
 A_j=\{\xi\in I_j:h_R(I_j)=M_{N,R}(\xi)\},\qquad
 E_{I_j}=A_j\setminus\bigcup_{i<j}A_i.
\]
These measurable sets form a partition of $\T$, because at each point at least one containing arc achieves the maximal value. Moreover,
\[
 E_{I_j}\subseteq I_j,\qquad
 M_{N,R}=h_R(I_j)\quad\text{on }E_{I_j}.
\]
Set $d_I=m(E_I)$. For each $J\in\mathcal D_N$,
\[
 \sum_{\substack{I\in\mathcal D_N\\I\subseteq J}}d_I
 =m\!\left(\bigcup_{\substack{I\in\mathcal D_N\\I\subseteq J}}E_I\right)
 \leq |J|.
\]
Hence $(d_I)$ is admissible, and
\[
 \sum_{I\in\mathcal D_N}d_Ih_R(I)
 =\sum_{I\in\mathcal D_N}\int_{E_I}M_{N,R}\,dm
 =\int_\T M_{N,R}\,dm.
\]
This proves equality \eqref{eq:endpointfinite}, and shows that the supremum is attained for the truncated finite problem.

Let $\mathcal A$ denote the collection of admissible sequences appearing in \eqref{eq:packingmax}. The restriction of any element of $\mathcal A$ to $\mathcal D_N$ satisfies the constraints defining $P_{N,R}$. Conversely, any vector admissible for $P_{N,R}$ extends to an element of $\mathcal A$ by assigning $d_I=0$ for indices outside $\mathcal D_N$. The constraints over $J\in\mathcal D_N$ remain unchanged; if $J$ has generation exceeding $N$, no arc from $\mathcal D_N$ is contained in $J$, so the associated sum vanishes. Therefore
\begin{equation}\label{eq:endpointextend}
 P_{N,R}=\sup_{d\in\mathcal A}
                     \sum_{I\in\mathcal D_N}d_Ih_R(I).
\end{equation}
All summands are nonnegative. Taking $R$ to be a positive integer and interchanging suprema,
\begin{align*}
 \sup_{d\in\mathcal A}\sum_{I\in\mathcal D}d_Ih(I)
 &=\sup_{d\in\mathcal A}\sup_{N,R}
                    \sum_{I\in\mathcal D_N}d_Ih_R(I)\\
 &=\sup_{N,R}P_{N,R}
 =\sup_{N,R}\int_\T M_{N,R}\,dm.
\end{align*}
The sequence $M_{n,n}$ increases pointwise to
$\sup_{I\ni\xi}h(I)$, and for any choice of $N,R$ there exists an integer
$n\geq\max\{N,R\}$ such that $M_{N,R}\leq M_{n,n}$. The monotone convergence theorem gives
\[
 \sup_{N,R}\int_\T M_{N,R}\,dm
 =\lim_{n\to\infty}\int_\T M_{n,n}\,dm
 =\int_\T\sup_{I\ni\xi}h(I)\,dm,
\]
which proves \eqref{eq:packingmax} without additional finiteness assumptions.
\end{proof}

We apply Lemma~\ref{lem:endpointpacking} with
$h(I)=\nu(Q_I^{\mathrm{up}})/|I|^3$. Using
\eqref{eq:treeadmissible} and \eqref{eq:hmeasure},
\begin{equation}\label{eq:endpointidentity}
 \|h_{\nu,2}^1\|_{B_1^2(\mathcal D)}
 =\int_\T\sup_{\substack{I\in\mathcal D\\\xi\in I}}
                  \frac{\nu(Q_I^{\mathrm{up}})}{|I|^3}\,dm(\xi).
\end{equation}
Substitute $\nu=\ind_E\mu_a^1$ and then take the supremum over $a\in\D$:
\begin{align*}
 \mathcal C_{1}^{1,p_2}(E)
 &=\sup_{a\in\D}\|h_{\ind_E\mu_a^1,2}^1\|_{B_1^2(\mathcal D)}\\
 &=\sup_{a\in\D}\int_\T
       \sup_{\substack{I\in\mathcal D\\\xi\in I}}
       \frac{\mu_a^1(E\cap Q_I^{\mathrm{up}})}{|I|^3}\,dm(\xi)
 =\mathcal M_{p_2}(E).
\end{align*}
This proves the identity in Theorem~\ref{thm:endpoint}, valid both for $E=\D$ and for $E=E_\rho$.

\section{Proofs of the main results}\label{sec:separation}
We first prove the separation theorem, then translate
boundedness and compactness of $\mathcal U$ and $\mathcal V$ into
capacity conditions, and finally deduce
Theorems~\ref{thm:subcritical}--\ref{thm:bloch}.

\subsection{Separation of the two terms}
Throughout this section $p_1>0$ and $0<p_2<1$ are fixed. We apply
the lemmas of Section~\ref{sec:prelim} with $p=p_1$ and choose $b$ and
$M$ accordingly. The following theorem separates the two terms.

\begin{theorem}\label{thm:separation}
Let $p_1>0$ and $0<p_2<1$, let $A,B$ be holomorphic functions on $\D$, and let
$\varphi:\D\to\D$ be holomorphic. Define
\[
 Tf=A(f\circ\varphi)+B(f'\circ\varphi),\quad
 T_0f=A(f\circ\varphi),\quad T_1f=B(f'\circ\varphi).
\]
Then $T:Q_{p_1}\to \mathcal N_{p_2}$ is bounded if and only if both $T_0$ and $T_1$ are bounded. In this case
\begin{equation}\label{eq:separationnorm}
 \|T_0\|+\|T_1\|\leq C_{p_1}\|T\|.
\end{equation}
Moreover, $T$ is compact if and only if both $T_0$ and $T_1$ are compact.
Consequently $W_{u,\varphi}:Q_{p_1}\to Q_{p_2}$ is bounded
(respectively compact) if and only if the two maps in \eqref{eq:product}
are bounded (respectively compact) as operators into $\mathcal N_{p_2}$.
\end{theorem}

We use the following compactness test for the operators in this theorem and for $W_{u,\varphi}$ \cite{WZ,Zhu}:
\begin{equation}\label{eq:localcriterion}
 \sup_n\|f_n\|_{Q_{p_1}}<\infty,\quad
 f_n\to0\text{ uniformly on compact subsets of }\D
 \quad\Longrightarrow\quad\|Sf_n\|\to0.
\end{equation}

For bounded $S$, \eqref{eq:localcriterion} is equivalent to compactness.
Indeed, norm convergence in $Q_{p_2}$ or $\mathcal N_{p_2}$ implies
locally uniform convergence, and the operators considered here preserve
locally uniform convergence. Conversely, by Montel's theorem and Fatou's
lemma, every bounded sequence in $Q_{p_1}$ has a subsequence converging
locally uniformly to some $f\in Q_{p_1}$, and \eqref{eq:localcriterion}
applies to the differences $f_n-f$.

\begin{proof}
Assume first that $T:Q_{p_1}\to \mathcal N_{p_2}$ is bounded. For each $a\in\D$,
introduce the positive Borel measure
\[
 \nu_a(E)=\int_{\varphi^{-1}(E)}|B(z)|^2
                  (1-|\varphi_a(z)|^2)^{p_2}\,dA(z).
\]
Fix a finite index set $\Lambda$ and coefficients
$c=(c_\lambda)_{\lambda\in\Lambda}$. Give $\{-1,1\}^\Lambda$ the uniform probability measure, and write $\varepsilon_\lambda$ for its coordinate random signs. Then
\[
 \mathbb E X=2^{-|\Lambda|}\sum_{\varepsilon\in\{-1,1\}^\Lambda}X(\varepsilon),
 \qquad
 \mathbb E(\varepsilon_\lambda\varepsilon_\eta)
 =\begin{cases}1,&\lambda=\eta,\\0,&\lambda\ne\eta.\end{cases}
\]
For $q=M,M+1$, set $f_{q,\varepsilon}=S_q(\varepsilon c)$. Sign changes do not alter the norm of the coefficient sequence, and Lemma~\ref{lem:atoms} gives
\[
 \sup_{q,\varepsilon}\|f_{q,\varepsilon}\|_{Q_{p_1}}
       \leq C_{p_1}\|c\|_{\mathfrak t_{p_1}}.
\]
Expanding the squared modulus and taking the expectation gives
\begin{align}\label{eq:randomidentity}
 \mathbb E|T f_{q,\varepsilon}(z)|^2
 &=\sum_{\lambda\in\Lambda}|c_\lambda|^2
   |A(z)h_{\lambda,q}(\varphi(z))
       +B(z)h'_{\lambda,q}(\varphi(z))|^2.
\end{align}
For $z\in\varphi^{-1}(V_\lambda)$, apply
Lemma~\ref{lem:jetmatrix} with $w=\varphi(z)$, $x=A(z)$, and $y=B(z)$.
Multiply through by $|c_\lambda|^2$ and integrate to obtain
\begin{align}
 &\sum_{\lambda\in\Lambda}|c_\lambda|^2
   \int_{\varphi^{-1}(V_\lambda)}
   \left(|A(z)|^2+\frac{|B(z)|^2}{(1-|a_\lambda|^2)^2}\right)
             (1-|\varphi_a(z)|^2)^{p_2}\,dA(z)\nonumber\\
 &\quad\leq C_M\sum_{q=M}^{M+1}\sum_{\lambda\in\Lambda}
       |c_\lambda|^2\int_\D|T h_{\lambda,q}(z)|^2
             (1-|\varphi_a(z)|^2)^{p_2}\,dA(z)\nonumber\\
 &\quad=C_M\sum_{q=M}^{M+1}\mathbb E\int_\D
       |T f_{q,\varepsilon}(z)|^2
             (1-|\varphi_a(z)|^2)^{p_2}\,dA(z)\nonumber\\
 &\quad\leq C_M\|T\|^2\sum_{q=M}^{M+1}
                   \mathbb E\|f_{q,\varepsilon}\|_{Q_{p_1}}^2
 \leq C_{p_1}\|T\|^2\|c\|_{\mathfrak t_{p_1}}^2.
 \label{eq:randomseparate}
\end{align}
Since $|A|^2\geq0$, we obtain
\begin{equation}\label{eq:diagonaltracebound}
 \sup_{a\in\D}\sum_{\lambda\in\Lambda}
 \frac{\nu_a(V_\lambda)}{(1-|a_\lambda|^2)^2}|c_\lambda|^2
 \leq C_{p_1}\|T\|^2\|c\|_{\mathfrak t_{p_1}}^2.
\end{equation}

Now take arbitrary $f\in Q_{p_1}$ and define
\[
 c_\lambda(f)=(1-|a_\lambda|^2)\sup_{w\in V_\lambda}|f'(w)|.
\]
By Lemma~\ref{lem:sampling} and monotonicity of the coefficient norm,
\[
 \|(c_{j,l}(f)\ind_{\{j\leq N\}})\|_{\mathfrak t_{p_1}}
 \leq\|c(f)\|_{\mathfrak t_{p_1}}
 \leq C_{p_1}\|f\|_{Q_{p_1},*}.
\]
Using \eqref{eq:diagonaltracebound} together with monotone convergence, for every $a\in\D$,
\[
 \sum_\lambda\frac{\nu_a(V_\lambda)|c_\lambda(f)|^2}
                         {(1-|a_\lambda|^2)^2}
 =\lim_{N\to\infty}\sum_{j\leq N,l}
       \frac{\nu_a(V_{j,l})|c_{j,l}(f)|^2}{(1-|a_{j,l}|^2)^2}
 \leq C_{p_1}\|T\|^2\|f\|_{Q_{p_1},*}^2.
\]
Hence
\begin{align*}
 \|T_1f\|_{\mathcal N_{p_2}}^2
 &=\sup_a\int_\D|f'(w)|^2\,d\nu_a(w)\\
 &\leq\sup_a\sum_\lambda
      \nu_a(V_\lambda)\sup_{V_\lambda}|f'|^2\\
 &=\sup_a\sum_\lambda
      \frac{\nu_a(V_\lambda)}{(1-|a_\lambda|^2)^2}|c_\lambda(f)|^2
 \leq C_{p_1}\|T\|^2\|f\|_{Q_{p_1},*}^2.
\end{align*}
It follows that
\[
 \|T_1\|\leq C_{p_1}\|T\|,\qquad
 \|T_0\|=\|T-T_1\|\leq(1+C_{p_1})\|T\|.
\]
The converse direction follows directly from the decomposition $T=T_0+T_1$.

Suppose now that $T$ is compact. Then both
$T_0$ and $T_1$ are bounded. Set
\begin{equation}\label{eq:epsilonJ}
 \varepsilon_J=\sup\left\{\|TS_qc\|_{\mathcal N_{p_2}}:
 \begin{array}{l}
 q\in\{M,M+1\},\ c\text{ has finite support},\\
 \|c\|_{\mathfrak t_{p_1}}\leq1,\quad\supp c\subseteq\{j\geq J\}
 \end{array}\right\}.
\end{equation}
Take any sequence $J_n\to\infty$, $q_n\in\{M,M+1\}$, and finitely supported coefficients $c^{(n)}$ satisfying
$\|c^{(n)}\|_{\mathfrak t_{p_1}}\leq1$ and
$\supp c^{(n)}\subseteq\{j\geq J_n\}$. Lemma~\ref{lem:atoms} gives
\[
 \sup_n\|S_{q_n}c^{(n)}\|_{Q_{p_1}}\leq C_{p_1},\qquad
 \sup_{|w|\leq R}|S_{q_n}c^{(n)}(w)|
 \leq C_{p_1,M,R}2^{-J_n(M-1)}\longrightarrow0\quad(R<1).
\]
Compactness combined with criterion \eqref{eq:localcriterion} gives
\[
 \|TS_{q_n}c^{(n)}\|_{\mathcal N_{p_2}}\longrightarrow0,
 \qquad\text{hence}\qquad\varepsilon_J\longrightarrow0.
\]

For finitely supported $c$ supported in $j\geq J$, homogeneity gives
\[
 \|TS_q(\varepsilon c)\|_{\mathcal N_{p_2}}
 \leq\varepsilon_J\|c\|_{\mathfrak t_{p_1}},\qquad q=M,M+1.
\]
Repeating the estimate \eqref{eq:randomseparate} with this bound gives
\begin{equation}\label{eq:compactdiagonal}
 \sup_a\sum_{j\geq J,\,l}
 \frac{\nu_a(V_{j,l})}{(1-|a_{j,l}|^2)^2}|c_{j,l}|^2
 \leq C_{p_1}\varepsilon_J^2\|c\|_{\mathfrak t_{p_1}}^2.
\end{equation}
Applying the bounded operator $T_1$ to the identity map $\id(w)=w$ gives
\begin{equation}\label{eq:uniformnumass}
 \sup_a\nu_a(\D)=\|T_1(\id)\|_{\mathcal N_{p_2}}^2<\infty.
\end{equation}

Let $(f_n)$ be a bounded sequence in $Q_{p_1}$ converging uniformly to zero on compact subsets. Write $r_J=1-2^{-J}$. The sets indexed by $j\geq J$ cover the region $\{|w|\geq r_J\}$. By Lemma~\ref{lem:sampling}, inequality \eqref{eq:compactdiagonal}, and monotone convergence,
\begin{align*}
 \|T_1f_n\|_{\mathcal N_{p_2}}^2
 &\leq \sup_{|w|\leq r_J}|f_n'(w)|^2\sup_a\nu_a(\D)
       +C_{p_1}\varepsilon_J^2\|f_n\|_{Q_{p_1},*}^2.
\end{align*}
For fixed $J$, Cauchy's integral formula guarantees
$\sup_{|w|\leq r_J}|f_n'(w)|\to0$. Therefore
\[
 0\leq\limsup_{n\to\infty}\|T_1f_n\|_{\mathcal N_{p_2}}^2
 \leq C_{p_1}\varepsilon_J^2\sup_n\|f_n\|_{Q_{p_1},*}^2
 \xrightarrow[J\to\infty]{}0.
\]
By \eqref{eq:localcriterion}, $T_1$ and $T_0=T-T_1$ are compact. Conversely, the sum of two compact operators is compact.

Finally set $A=u'$ and $B=u\varphi'$. Differentiation carries a bounded or compact $W_{u,\varphi}$ to an operator $T$ with the same property. Conversely,
\begin{equation}\label{eq:volterradecomposition}
 W_{u,\varphi}f(z)=u(0)f(\varphi(0))
       +\int_0^z\mathcal U f(\xi)\,d\xi
       +\int_0^z\mathcal V f(\xi)\,d\xi.
\end{equation}
Integration defines a bounded mapping from $\mathcal N_{p_2}$ into $\mathring Q_{p_2}$, while the initial term is a rank-one bounded operator. This completes the proof of the final claim.
\end{proof}

\subsection{Compactness of the two embeddings}
By Theorem~\ref{thm:separation},
\[
 W_{u,\varphi}\text{ compact}\quad\Longleftrightarrow\quad
 \mathcal U,\mathcal V:Q_{p_1}\to\mathcal N_{p_2}\text{ compact}.
\]
We now express these conditions through the two capacities.

Substitute $p=p_1$ and $\nu=\mu_a^1$ into \eqref{eq:localtrace} and take the supremum over all $a$. This gives
\begin{equation}\label{eq:derivativeembedding}
 \|\mathcal V:Q_{p_1}\to\mathcal N_{p_2}\|^2
 \asymp\mathcal C_{1}^{p_1,p_2}(\D).
\end{equation}

\begin{proposition}\label{prop:derivcompact}
Assume that $\mathcal V:Q_{p_1}\to \mathcal N_{p_2}$ is bounded. Then
\begin{equation}\label{eq:derivativecompacttail}
 \mathcal V\text{ is compact}
 \quad\Longleftrightarrow\quad
 \lim_{r\uparrow1}\sup_{a\in\D}\bigl\|h_{\ind_{\{|w|>r\}}\mu^1_a,2}^{p_1}\bigr\|_{B_{p_1}^2(\mathcal D)}=0.
\end{equation}
\end{proposition}
\begin{proof}
Suppose $\mathcal V$ is compact. Define
\[
 \varepsilon_J=\sup\bigl\{\|\mathcal V S_Mc\|_{\mathcal N_{p_2}}:
 c\text{ finite},\ \|c\|_{\mathfrak t_{p_1}}\leq1,\ \supp c\subseteq\{j\geq J\}\bigr\}.
\]
Combining \eqref{eq:atomsbound}, \eqref{eq:atomstail}, and \eqref{eq:localcriterion}, we obtain $\varepsilon_J\to0$ as $J\to\infty$. Consequently, for every finitely supported coefficient sequence $c$ supported in $j\geq J$,
\[
 \|\mathcal V S_Mc\|_{\mathcal N_{p_2}}\leq\varepsilon_J\|c\|_{\mathfrak t_{p_1}}.
\]
For coefficients supported on indices $j\geq J$, averaging over random sign choices together with the pointwise lower bound $|h'_{\lambda,M}|\geq c_M/(1-|a_\lambda|^2)$ on $V_\lambda$ gives
\begin{align*}
 \sum_{j\geq J,l}\frac{\mu^1_a(V_{j,l})}{(1-|a_{j,l}|^2)^2}|c_{j,l}|^2
 &\leq C_M\mathbb E\int_\D
       \bigl|u(z)\varphi'(z)\bigl(S_M(\varepsilon c)\bigr)'(\varphi(z))\bigr|^2
       (1-|\varphi_a(z)|^2)^{p_2}\,dA(z)\\
 &\leq C_M\varepsilon_J^2\|c\|_{\mathfrak t_{p_1}}^2.
\end{align*}
The union of these sets is $\{|w|\geq r_J\}$, where $r_J=1-2^{-J}$. The
definition \eqref{eq:coefficientgauge} and \eqref{eq:treecomparison}
therefore yield
\[
 \sup_a\bigl\|h_{\ind_{\{|w|\geq r_J\}}\mu^1_a,2}^{p_1}\bigr\|_{B_{p_1}^2(\mathcal D)}
 \leq C_{p_1}\varepsilon_J^2\longrightarrow0.
\]
Monotonicity of the expression in $r$ extends this limit to any sequence $r\uparrow1$.

Conversely, assume the limit in \eqref{eq:derivativecompacttail} holds. Testing the bounded operator $\mathcal V$ against the function $f(w)=w$ gives
\[
 M_1:=\sup_a\mu^1_a(\D)=\|u\varphi'\|_{\mathcal N_{p_2}}^2<\infty.
\]
Let $\{f_n\}$ be a bounded sequence in $Q_{p_1}$ that converges to zero uniformly on every compact subset of $\D$. For fixed $r<1$, the comparison \eqref{eq:localtrace} gives
\begin{align*}
 \|\mathcal V f_n\|_{\mathcal N_{p_2}}^2
 &=\sup_a\int_\D |f_n'(w)|^2\,d\mu^1_a(w)\\
 &\leq M_1\sup_{|w|\leq r}|f_n'(w)|^2
 +C_{p_1}\|f_n\|_{Q_{p_1},*}^2
        \sup_a\bigl\|h_{\ind_{\{|w|>r\}}\mu^1_a,2}^{p_1}\bigr\|_{B_{p_1}^2(\mathcal D)}.
\end{align*}
Writing $C=\sup_n\|f_n\|_{Q_{p_1},*}^2<\infty$, we obtain
\[
 \limsup_{n\to\infty}\|\mathcal V f_n\|_{\mathcal N_{p_2}}^2
 \leq C_{p_1} C\sup_a\bigl\|h_{\ind_{\{|w|>r\}}\mu_a^1,2}^{p_1}\bigr\|_{B_{p_1}^2(\mathcal D)}
 \xrightarrow[r\uparrow1]{}0.
\]
Apply the test \eqref{eq:localcriterion}.
\end{proof}

For any Borel subset $E\subseteq\D$, define
\begin{equation}\label{eq:L0}
 \mathcal L_EF(z)=u'(z)\int_E
       L_b(\varphi(z),w)F(w)\,dA_b(w).
\end{equation}
We substitute $p=p_1$ and $\nu=\mu^0_{a,t}$ into \eqref{eq:restrictedembedding} and then take the supremum over the parameters $a$ and $t$. This gives
\begin{equation}\label{eq:transfer}
 \|\mathcal L_E:T_{p_1}\to\mathcal N_{p_2}\|^2
 \asymp\mathcal C_{0}^{p_1,p_2}(E).
\end{equation}
The implicit constants are independent of $E$, $a$, and $t$. In the special case $E=\D$, the reconstruction formula together with the value embedding theorem further gives
\begin{equation}\label{eq:valueembedding}
 \|\mathcal U:\mathring Q_{p_1}\to\mathcal N_{p_2}\|^2
 \asymp\mathcal C_{0}^{p_1,p_2}(\D).
\end{equation}

\begin{proposition}\label{prop:valuecompact}
Suppose $\mathcal U:Q_{p_1}\to \mathcal N_{p_2}$ is bounded. Then
\begin{equation}\label{eq:valuecompact}
 \mathcal U\text{ is compact}\quad\Longleftrightarrow\quad
                   \mathcal C_{0}^{p_1,p_2}(E_\rho)\longrightarrow0\quad(\rho\uparrow1).
\end{equation}
\end{proposition}
\begin{proof}
First assume that $\mathcal U$ is compact. Lemma~\ref{lem:taillocal} supplies the estimates
\begin{gather*}
 \sup_{\rho<1}\sup_{\|F\|_{T_{p_1}}\leq1}\|P_\rho F\|_{Q_{p_1}}<\infty,\\
 \sup_{\|F\|_{T_{p_1}}\leq1}\sup_{|\zeta|\leq R}|P_\rho F(\zeta)|
 \leq C_{p_1,b,R}\eta_b(\rho)\longrightarrow0\qquad(R<1).
\end{gather*}
If $\sup_{\|F\|_{T_{p_1}}\leq1}\|\mathcal U P_\rho F\|_{\mathcal N_{p_2}}$ failed to vanish as $\rho\uparrow1$, we could select
\[
 \rho_n\uparrow1,\quad \|F_n\|_{T_{p_1}}\leq1,\quad
 \|\mathcal U P_{\rho_n}F_n\|_{\mathcal N_{p_2}}\geq\varepsilon>0
\]
for some fixed $\varepsilon>0$. However, the preceding bounds combined with \eqref{eq:localcriterion} force $\|\mathcal U P_{\rho_n}F_n\|_{\mathcal N_{p_2}}\to0$, which gives a contradiction. It follows that
\[
 \mathcal C_{0}^{p_1,p_2}(E_\rho)
 \asymp\|\mathcal L_{E_\rho}\|_{T_{p_1}\to \mathcal N_{p_2}}^2
 =\|\mathcal U P_\rho\|_{T_{p_1}\to \mathcal N_{p_2}}^2\longrightarrow0.
\]

For the converse direction, Lemma~\ref{lem:interiorcompact} guarantees that $S_\rho$ is compact. Moreover,
\begin{align*}
 (\mathcal U-\mathcal U S_\rho)f
 &=\mathcal U P_\rho(f')=\mathcal L_{E_\rho}(f'),\\
 \|\mathcal U-\mathcal U S_\rho\|_{Q_{p_1}\to \mathcal N_{p_2}}^2
 &\leq C_{p_1}\|\mathcal L_{E_\rho}\|_{T_{p_1}\to \mathcal N_{p_2}}^2
 \leq C_{p_1}\mathcal C_{0}^{p_1,p_2}(E_\rho)\longrightarrow0.
\end{align*}
Thus $\mathcal U S_\rho\to\mathcal U$ in operator norm. Since each
$\mathcal U S_\rho$ is compact, $\mathcal U$ is compact.
\end{proof}

\subsection{The capacity criterion and the main theorems}
We first combine the two embedding conditions. The endpoint identity
and the comparison below then give the three main theorems.

\begin{proposition}\label{prop:unified}
Let $p_1>0$ and $0<p_2<1$. Suppose $u\in\Hol(\D)$ and $\varphi:\D\to\D$ is holomorphic.
Then $W_{u,\varphi}:Q_{p_1}\to Q_{p_2}$ is bounded if and only if
\begin{equation}\label{eq:boundediff}
 u\in Q_{p_2},\qquad \mathcal C_{0}^{p_1,p_2}(\D)+\mathcal C_{1}^{p_1,p_2}(\D)<\infty.
\end{equation}
In that case,
\begin{align}\label{eq:mainnorm}
 \|W_{u,\varphi}\|^2\asymp{}&
 |u(0)|^2\log^2\frac{e}{1-|\varphi(0)|^2}+\|u\|_{Q_{p_2},*}^2\nonumber\\
 &+\mathcal C_{0}^{p_1,p_2}(\D)+\mathcal C_{1}^{p_1,p_2}(\D).
\end{align}
A bounded $W_{u,\varphi}$ is compact if and only if
\begin{equation}\label{eq:compactiff}
 \lim_{\rho\uparrow1}\bigl[\mathcal C_{0}^{p_1,p_2}(E_\rho)+\mathcal C_{1}^{p_1,p_2}(E_\rho)\bigr]=0.
\end{equation}
All comparison constants depend only on $p_1,p_2$ and the fixed reconstruction parameter.
\end{proposition}
\begin{proof}
Write $W=W_{u,\varphi}$. If $W$ is bounded, separation and the two
embedding comparisons yield
\begin{align*}
 u=W1&\in Q_{p_2},\qquad \|u\|_{Q_{p_2},*}\leq\|W\|,\\
 \mathcal C_{0}^{p_1,p_2}(\D)+\mathcal C_{1}^{p_1,p_2}(\D)
 &\lesssim\|\mathcal U\|^2+\|\mathcal V\|^2\lesssim\|W\|^2,\\
 |u(0)|\log\frac{e}{1-|\varphi(0)|^2}
 &\asymp |u(0)|\|\operatorname{ev}_{\varphi(0)}\|_{Q_{p_1}^*}
 =\|\operatorname{ev}_0 W\|\leq\|W\|.
\end{align*}
These are the lower estimates in \eqref{eq:mainnorm}.

Conversely, suppose condition \eqref{eq:boundediff} holds. Using the decomposition $f=f(0)+P_{b,0}(f')$, we obtain
\begin{align*}
 \|\mathcal U f\|_{\mathcal N_{p_2}}
 &\leq |f(0)|\|u\|_{Q_{p_2},*}
       +C_{p_1}\mathcal C_{0}^{p_1,p_2}(\D)^{1/2}\|f\|_{Q_{p_1},*},\\
 \|\mathcal V f\|_{\mathcal N_{p_2}}
 &\leq C_{p_1}\mathcal C_{1}^{p_1,p_2}(\D)^{1/2}\|f\|_{Q_{p_1},*},\\
 |(Wf)(0)|
 &\leq C_{p_1}|u(0)|\log\frac{e}{1-|\varphi(0)|^2}\|f\|_{Q_{p_1}}.
\end{align*}
Therefore
\begin{align*}
 \|Wf\|_{Q_{p_2}}
 &\leq |(Wf)(0)|+\|\mathcal U f\|_{\mathcal N_{p_2}}
                            +\|\mathcal V f\|_{\mathcal N_{p_2}}\\
 &\leq C_{p_1,p_2}\left[
       |u(0)|\log\frac{e}{1-|\varphi(0)|^2}+\|u\|_{Q_{p_2},*}
       +\sum_{j=0}^1\mathcal C_{j}^{p_1,p_2}(\D)^{1/2}\right]\|f\|_{Q_{p_1}}.
\end{align*}
Square this estimate and take the supremum over the unit ball.
Assuming boundedness, Theorem~\ref{thm:separation} and Propositions~\ref{prop:derivcompact}--\ref{prop:valuecompact} imply
\begin{align*}
 W\text{ compact}
 &\iff \mathcal U,\mathcal V:Q_{p_1}\to \mathcal N_{p_2}\text{ are compact}\\
 &\iff \mathcal C_{0}^{p_1,p_2}(E_\rho)\to0
                  \quad\text{and}\quad\mathcal C_{1}^{p_1,p_2}(E_\rho)\to0\\
 &\iff\lim_{\rho\uparrow1}
       \bigl[\mathcal C_{0}^{p_1,p_2}(E_\rho)+\mathcal C_{1}^{p_1,p_2}(E_\rho)\bigr]=0.
\end{align*}
The last equivalence uses nonnegativity of both quantities.
\end{proof}

Both limits are taken after the suprema in \eqref{eq:Cglobal}.

\begin{corollary}\label{cor:blochderivative}
If $p_1>1$, then for every Borel subset $E\subseteq\D$,
\[
 \mathcal C_{1}^{p_1,p_2}(E)\asymp
 \sup_{a\in\D}\int_{\varphi^{-1}(E)}
 \frac{|u(z)\varphi'(z)|^2}{(1-|\varphi(z)|^2)^2}
             (1-|\varphi_a(z)|^2)^{p_2}\,dA(z).
\]
This applies both to $E=\D$ and to $E=E_\rho$.
\end{corollary}
\begin{proof}
Let $p_1>1$. Every admissible sequence satisfies $d_I\leq|I|^{p_1}$.
For the reverse estimate, set $d_I=(1-2^{1-p_1})|I|^{p_1}$. For any interval $J$,
\begin{align*}
 |J|^{-p_1}\sum_{I\subseteq J}d_I
 &=(1-2^{1-p_1})|J|^{-p_1}
   \sum_{n=0}^\infty\sum_{\substack{I\subseteq J\\|I|=2^{-n}|J|}}|I|^{p_1}\\
 &=(1-2^{1-p_1})\sum_{n=0}^\infty2^n2^{-np_1}=1.
\end{align*}
We therefore have, with values taken in $[0,\infty]$,
\begin{equation}\label{eq:blochcapacitysum}
 (1-2^{1-p_1})\sum_I\frac{\nu(Q_I^{\mathrm{up}})}{|I|^2}
 \leq\|h_{\nu,2}^{p_1}\|_{B_{p_1}^2(\mathcal D)}
 \leq\sum_I\frac{\nu(Q_I^{\mathrm{up}})}{|I|^2}.
\end{equation}
For any point $w\in Q_I^{\mathrm{up}}$,
\[
 \frac{|I|}{2}<1-|w|\leq|I|,\qquad
 \frac{|I|}{2}<1-|w|^2\leq2|I|.
\]
The upper-half tiles $Q_I^{\mathrm{up}}$ form a partition of the unit disk. It follows that
\begin{align*}
 \frac14\int_\D\frac{d\nu(w)}{(1-|w|^2)^2}
 &\leq\sum_I\frac{\nu(Q_I^{\mathrm{up}})}{|I|^2}
 \leq4\int_\D\frac{d\nu(w)}{(1-|w|^2)^2}.
\end{align*}
Now apply this estimate to the measure $\nu=\ind_E\mu_a^1$ and invoke the pullback identity:
\begin{align*}
 \int_\D\frac{d(\ind_E\mu_a^1)(w)}{(1-|w|^2)^2}
 &=\int_{\varphi^{-1}(E)}
       \frac{|u(z)\varphi'(z)|^2}{(1-|\varphi(z)|^2)^2}
                 (1-|\varphi_a(z)|^2)^{p_2}\,dA(z).
\end{align*}
Taking the supremum over all $a\in\D$, we obtain, uniformly with respect to $E$,
\begin{equation}\label{eq:blochcapacityconstants}
 \frac{1-2^{1-p_1}}4\,\mathcal J_{p_2}(E)
 \leq\mathcal C_{1}^{p_1,p_2}(E)\leq4\mathcal J_{p_2}(E).
\end{equation}
\end{proof}

\begin{proof}[Proof of Theorems~\ref{thm:subcritical}--\ref{thm:bloch}]
When $0<p_1<1$, Proposition~\ref{prop:unified} directly gives Theorem~\ref{thm:subcritical}. For $p_1=1$, Lemma~\ref{lem:endpointpacking} together with \eqref{eq:endpointidentity} gives
\[
 \mathcal C_{1}^{1,p_2}(\D)=\mathcal M_{p_2}(\D),\qquad
 \mathcal C_{1}^{1,p_2}(E_\rho)=\mathcal M_{p_2}(E_\rho).
\]
Substitution into Proposition~\ref{prop:unified} proves Theorem~\ref{thm:endpoint}, for both boundedness and compactness.
If $p_1>1$, the identification
\[
 Q_{p_1}=Q_2=\mathcal B,\qquad
 \|f\|_{Q_{p_1}}\asymp_{p_1}\|f\|_{Q_2},
\]
shows that the boundedness and compactness properties remain unchanged upon replacing $Q_{p_1}$ by $Q_2$.
Applying Proposition~\ref{prop:unified} with $p_1=2$ and using \eqref{eq:blochcapacityconstants}, we obtain
\begin{gather*}
 W\text{ is bounded}\iff u\in Q_{p_2},\quad
             \mathcal C_{0}^{2,p_2}(\D)+\mathcal J_{p_2}(\D)<\infty,\\
 W\text{ is compact}\iff
 \lim_{\rho\uparrow1}\bigl[\mathcal C_{0}^{2,p_2}(E_\rho)
                              +\mathcal J_{p_2}(E_\rho)\bigr]=0
 \quad\text{provided }W\text{ is bounded}.
\end{gather*}
In this argument, the value capacity is defined with $p=2$ and with
$b$ chosen accordingly.
\end{proof}

\section*{Declarations}

\noindent\textbf{Funding.} This work was supported by the National
Natural Science Foundation of China (Grant No.~12371131).

\medskip

\noindent\textbf{Conflict of interest.} The authors declare no
competing interests.

\medskip

\noindent\textbf{Data availability.} Data sharing is not applicable
for this article as no datasets were generated or analyzed during the
current study.

\medskip

\noindent\textbf{Use of Generative-AI tools declaration.} During the
preparation of this work, the authors used an artificial-intelligence
tool for language editing, \LaTeX{} formatting and local mathematical
reasoning. All mathematical content was written, checked, and approved
by the authors, who take full responsibility for the content of the
paper.

\end{document}